\documentclass[11pt]{article}
\usepackage{amsmath,amstext,amsthm,amsfonts,amscd,amssymb}
\usepackage[dvips]{graphicx}
\usepackage{epic,eepic}
\usepackage[latin1]{inputenc}
\theoremstyle{plain}
\newtheorem{theorem}{Theorem }[section]

\newtheorem{lemma}[theorem]{Lemma}

\newtheorem{maintheorem}{Theorem}

\theoremstyle{definition}
\newtheorem{remark}[theorem]{Remark}

\newtheorem{definition}[theorem]{Definition}

\newcommand{\field}[1]{\mathbb{#1}}
\newcommand{\real}{\field{R}}

\renewcommand{\natural}{\field{N}}

\newcommand{\al} {\alpha}       
\newcommand{\be} {\beta}        
\newcommand{\ga} {\gamma}    
\newcommand{\de} {\delta}       \newcommand{\De}{\Delta}
\newcommand{\ep} {\epsilon}     

\newcommand{\ze} {\zeta}
\newcommand{\te} {\theta}

\newcommand{\la} {\lambda}

\newcommand{\si} {\sigma}

\newcommand{\vfi}{\varphi}

\newcommand{\cqd}{\hfill q.e.d. \vspace{.1in} }

\begin{document}

\title{\Large{\textbf{ Immersions with fractal set of points of zero 
Gauss-Kronecker curvature }}}

\author{ Alexander Arbieto\footnote{A. Arbieto is supported by CNPq/Brazil}
\text{ and}
Carlos Matheus\footnote{C. Matheus is supported by Faperj/Brazil} }

\date{February 25, 2003}

\maketitle

\begin{abstract} We construct, for any ``good'' Cantor set $F$ of $S^{n-1}$, an
immersion of the sphere $S^n$ with set of points of zero Gauss-Kronecker 
curvature equal to $F\times D^{1}$, where $D^{1}$ is the $1$-dimensional
disk. In particular these examples show that the theorem of 
Matheus-Oliveira strictly extends two results by 
do Carmo-Elbert and 
Barbosa-Fukuoka-Mercuri.   
\end{abstract}

\section{Introduction}
Let $x:M^n\rightarrow\real^{n+1}$ be a codimension one Euclidean immersion of an
orientable manifold $M$. Let $N:M^n\rightarrow S^n$ be the associated Gauss map
in the given orientation. Recall that $A=dN$ is self-adjoint and its eigenvalues
are the principal curvatures. We denote $H_n=det(dN)$ the Gauss-Kronecker
curvature and $rank(x):=rank(N):=\min\limits_{p\in M}
rank(d_p N)$.

A compact set $F\subset S^n$ is called a \emph{good} Cantor set if 
$S^n-F=\bigcup_{i\in\natural}U_i$,
is the disjoint union of open balls $U_i$ in $S^n$ (with the standard metric) of
radius bounded by a small constant $\de_0=\de_0(n)$ (to be choosen later).

Our main result is : 

\begin{maintheorem}\label{imersao fractal} For any $F\subset S^{n-1}$
a good Cantor set, there are immersions $x:S^n\rightarrow\real^{n+1}$ such that
$rank(x)=n-1$, the Gauss-Kronecker curvature is non-negative and $\{p\in S^n: 
H_n(p)=0
\}=F\times D^{1}$, where $D^{1}$ is the $1$-dimensional disk.
\end{maintheorem}

Before starting some comments, we briefly recall the definition :

\begin{definition} A complete orientable hypersurface $M$ has finite geometrical
type (in~\cite{BFM} sense) if $M$ is equal to a compact manifold
$\overline{M}$ minus a finite number of points (called \emph{ends}), the
Gauss map $N:M\rightarrow S^n$ extends continuously to each ``end'' and the set
of points of zero Gauss-Kronecker curvature is contained in a finite union of
submanifolds of dimension $\leq n-2$ ; $M$ 
has finite total curvature if $\int_M |A|^n <\infty$, where $|A|
=( \sum\limits_{i=1}^{n}k_i^2 )^{1/2}$, $k_i$ are the principal curvatures. We
remark that, as showed in \cite{dCE}, if $M$ has finite total curvature then $M$
has finite geometrical type. 
\end{definition}

The motivation of our theorem~\ref{imersao fractal} are the following results. 
In a recent work, do Carmo and Elbert~\cite{dCE} show that :

\begin{theorem}[do Carmo, Elbert] If $M$ is a 
hypersurface of finite total curvature with Gauss-Kronecker curvature 
$H_n\neq 0$ everywhere then $M$ is topologically a sphere 
minus a finite number of points. 
\end{theorem}

In fact, this result can 
be
improved, as showed in~\cite{BFM} : 

\begin{theorem}[Barbosa, Fukuoka, Mercuri]If $M$ is a $2n$ dimensional 
hypersurface 
of finite geometrical type such that $\{ p\in M: H_n(p)=0 \}$ is a  
subset of a finite union of submanifolds with dimension less than $n-1$ then 
the hypersurface is a sphere minus \emph{two} points. If $M$ is minimal, $M$ is
the $2n$-catenoid.
\end{theorem} 

In general, it is not
easy to obtain the hypothesis of the result by Barbosa, Fukuoka and Mercuri for
arbitrary immersions, since the classical theorems (Sard, Moreira~\cite{M}) 
treat only
the critical values (in our case, the spherical image of points of zero
curvature).
 With this difficult in mind, Matheus and 
Oliveira~\cite{MO},
using the concept of Hausdorff dimension, generalize the previous results : 

\begin{theorem}[Matheus, Oliveira] If $x:M\rightarrow\real^{n+1}$ is an
immersion of finite geometrical type such that $rank(x)\geq k$ and the $(k-[\frac{n}{2}])$- dimensional 
Hausdorff measure of $N(\{ p\in M:
H_n(p)=0 \} )$ is zero then $M$ is topologically a sphere minus a finite number of
points ($[r]=\text{integer part of } r$). If $M$ is minimal, $M$ is the
catenoid.
\end{theorem}

However, it is not obvious the existence of immersions satisfying
Matheus-Oliveira's hypothesis which does not satisfies the assumptions of 
Barbosa-Fukuoka-Mercuri. The question about existence of ``fractal immersions'' 
was posed to second author by Walcy Santos during the
Differential Geometry seminar at IMPA. The main goal of this paper is to show
how one can construct immersions with ``large'' fractal set of points of zero
curvature. Clearly, the existence of such immersions follows from our
theorem~\ref{imersao fractal}. The outline of proof of theorem~\ref{imersao
fractal} is :

\begin{itemize}
\item First, we construct immersions $\vfi$ of $S^n$ with rank equal to $k$, 
$\{ H_n=0
\} := \{p
:H_n(p)=0 \}=S^k\times D^{n-k}$ ($D^{n-k}$ is a $(n-k)$-dimensional disk and
$S^k$ is a round sphere) and $N(\{ H_n=0 \})=S^k$; 

\item Second, if $U$ is a ball in $S^{n-1}$, we produce a modification
$\vfi|_U$ of $\vfi$ such that $\vfi_U=\vfi$ in $(S^{n-1}-U)\times
D^{1}$ and the Gauss-Kronecker curvature of $\vfi_U(U\times D^{1})$ does not 
vanishes;

\item Finally, we consider $F=S^{n-1}-\bigcup\limits_{i=1}^{\infty} U_i$ and we produce the 
desired immersion
by induction, using the two steps above.

\end{itemize}

In next three sections, we are going to make precise the steps described
above. In other words, we describe explicitly the immersions with the properties
commented in the previous outline.  

To finish this introduction, we observe that in the special case of ``good'' Cantor
sets, even if it has positive measure, the manifold can be the sphere (this
occurs because the ``singular'' set has ``good'' geometry). In particular, the
theorems of Matheus-Oliveira are not sharp. Moreover, the proof of our
theorem~\ref{imersao fractal} shows that the round balls can be replaced by sets
with a well-defined ``distance function''. However, the proof only works in
codimension 1 (i.e., for good Cantor sets $F\subset S^{n-1}$), by technical 
reasons. Also, we can apply these technics in the non-compact case. This is showed in last 
section of this paper.  

\section{Some immersions with cylindrical pieces}

The main result of this section is the following lemma :

\begin{lemma}\label{im. cilindrica} There are immersions 
$\vfi:S^n\rightarrow\real^{n+1}$ such that
$rank(\vfi)=k$ and $\{p\in S^n: H_n(p)=0\}=S^k\times D^{n-k}$, where $S^k$ is a
round sphere in $\real^{k+1}$ with radius $0<\sqrt{\ga}<1$, $D^{n-k}$ is the
$(n-k)$-dimensional disk of radius $\sqrt{\al}$, $\al < \frac{1}{2}$ and $\al < \ga < 1-\al$.  
\end{lemma}

The ideia of the proof of this lemma is flatten the boundary of a hemisphere 
such that the
curvature is positive everywhere except at the boundary. Now take a
cylinder and glue isometrically the boundaries of the cylinder and the
hemisphere (see figure 1 after the proof of lemma~\ref{im. cilindrica}). 
Gromov~\cite{G} uses this ideia (in other context) to construct some
examples of manifolds of nonpositive curvature with special properties. Because
the authors does not know any reference in literature where these examples are
constructed in details (the known references deals only with higher codimension
surgeries~\cite{GL}), we present the proof of the 
lemma~\ref{im. cilindrica}.

\proof[Proof of lemma~\ref{im. cilindrica}] Fix some $\al<\be<\ga$. We write 
$\real^{n+1}\ni p=(x,y)$, where
$x\in\real^{k+1}$ and $y\in\real^{n-k}$. We denote $||.||$ the Euclidean metric.
Consider some real function $\nu\in C^{\infty}$ s.t. $\nu(r)\equiv \ga$ if
$0\leq r\leq \al$, $\nu(r)\equiv 1$ if $\be\leq r\leq 1$ and $\nu$ is
strictly increasing in $[\al,\be]$. The immersion $\vfi:S^n\rightarrow
\real^{n+1}$ is:
$$\vfi(x,y)=(\te(||y||^2)\cdot x, y)=(z,w),\textrm{ where }
\te(r)=\sqrt{\frac{\nu(r)-r\cdot\mu(r)}{1-r}},$$
and $\mu$ is a convenient real
function (it will be defined later). We take $\mu$ such that $\mu(r)\equiv 0$ if
$0\leq r\leq\al$ and $\mu(r)\equiv 1$ if $\be\leq r\leq 1$. These implies that
$\vfi(x,y)=(x,y)$ if $||y||^2\geq\be$,
$\vfi(x,y)=(\frac{\sqrt{\ga}}{\sqrt{1-||y||^2}}\cdot x, y)$  if $||y||^2\leq\al$. 
In other words, $\vfi$ has
a spherical piece and a cylindrical piece. Now, it is sufficient to define 
$\vfi$ (i.e., $\mu$) in such
way that the Gauss-Kronecker curvature is positive except at the cylindrical
piece. 

By definition, using $y=w$, $(x,y)\in S^n\Rightarrow ||x||^2+||y||^2 =1$,
we have $||z||^2 + \mu(||w||^2)\cdot ||w||^2 = \nu(||w||^2)$.
That's it, if 
$$f(z,w)=||z||^2 + \mu(||w||^2)\cdot ||w||^2 - \nu(||w||^2),
$$ 
then
$f|_{\vfi(S^n)}\equiv 0$. As usual, the Gauss map is
$N(z,w)=\frac{\textrm{grad} \ f}{||\textrm{grad} \ f||}$. But $\frac{\partial f}{\partial z_i} =
2\cdot z_i$ and $\frac{\partial f}{\partial w_j} = 2\cdot w_j\cdot\{\mu +
||w||^2\cdot\mu'-\nu'\}$. For sake of simplicity, we will denote $c_1(r)=\mu + r\cdot\mu'-\nu'$.
Then, $||\textrm{grad} \ f||=2\cdot\sqrt{\nu + ||w||^2\cdot(c_1^2-\mu)}$. We denote
$c_2(r)=\sqrt{\nu + r\cdot(c_1^2-\mu)}$. Now we have: 
$$ N(z,w)=\frac{1}{c_2(||w||^2)}\cdot (z, \ c_1(||w||^2)\cdot w).
$$ 
Clearly, $\frac{\partial N_l}{\partial z_i}=0$ if $l\neq i$, $\frac{\partial
N_l}{\partial z_i}=\frac{1}{c_2}$ if $l=i$. Analogously, $\frac{\partial
N_l}{\partial w_j}=2\cdot(\frac{1}{c_2})'\cdot z_l\cdot w_j$ if $l\leq k+1$ and,
if $l\geq k+2$, 

\begin{displaymath}
\frac{\partial N_l}{\partial w_j}= \left\{ \begin{array}{ll}
2\cdot(\frac{c_1}{c_2})'\cdot w_j\cdot w_l & \textrm{if $l\neq j$} \\
2\cdot(\frac{c_1}{c_2})'\cdot w_j\cdot w_l + \frac{c_1}{c_2} & \textrm{if $l=j$}
\\
\end{array} \right.
\end{displaymath} 

These
computations implies that :

\begin{displaymath}
\mathbf{dN} =
\left( \begin{array}{cc}
\frac{1}{c_2} \cdot I_{k+1} & \bigstar \\
0 & A + \frac{c_1}{c_2}\cdot I_{n-k} \\
\end{array} \right),
\end{displaymath}
where $I_m$ is the $(m\times m)$-identity matrix and $A=2\cdot
(\frac{c_1}{c_2})'\cdot [w_j\cdot w_l]_{jl}$ ($w_j$ is the $j$-th component of 
$w$).
We observe that $A$ has eigenvalues $0$ with multiplicity
$(n-k-1)$ and $2\cdot (\frac{c_1}{c_2})'\cdot ||w||^2$ with multiplicity $1$. In 
particular, the principal curvatures of $\vfi$ are : $\frac{1}{c_2}$ with
multiplicity $k$, $\frac{c_1}{c_2}$ with multiplicity $n-k-1$ and $2\cdot
(\frac{c_1}{c_2})'\cdot ||w||^2 + \frac{c_1}{c_2}$ with multiplicity $1$.
Because 
$\det
(dN-\la\cdot Id) = \det (\frac{1}{c_2} \cdot Id) \cdot \det (A + \frac{c_1}{c_2}\cdot Id -
\la\cdot Id)$.

We need
$H_n>0$ if $||y||^2=||w||^2 >\al$, $H_n=0$ if $||w||^2\leq\al$. To solve this,
we consider $\psi\in C^{\infty}$, $\psi|_{[0,\al]}\equiv 0$ and $\psi|_{[\be,1]}
\equiv 1$. The admissible $\mu$ needs to satisfy:
\begin{eqnarray}
\frac{c_1}{c_2}=\rho \\ 
2 r\rho' + \rho = \psi
\end{eqnarray}
where $\rho|_{[0,\al]}\equiv 0$, $\rho|_{[\be,1]} \equiv
1$ also. As we will see
later, this property is sufficient to conclude the proof.

In order to make 
the following rigorous, fix any 
$0<\widetilde{\al}<\al$
and $\be<\widetilde{\be}<1$. The ODE (2) can 
be
explicitly solved : $(2)\iff \rho' + \frac{\rho}{2r} = \frac{\psi}{2r}\iff \big(
e^{\int_{\al}^{r}\frac{1}{2s} \ ds}\cdot\rho \big)'=\frac{\psi}{2r}$.
Integrating (and using $\rho(\al)=0$), we get $(2)\Leftrightarrow\rho(r)=
\frac{1}{\sqrt{2r}}\cdot\int\limits_{\al}^{r} \frac{\psi}{\sqrt{2s}} \ ds$, for
$r\in [\widetilde{\al}, \widetilde{\be}]$. 

Now, we solve the (implicit) ODE (1). Observe that 
$(1)\Leftrightarrow c_1^2=\rho^2\nu + c_1^2 r\rho - \mu
r\rho^2 \Leftrightarrow \frac{c_1^2}{\nu - r\mu} = \frac{\rho^2}{1 - r\rho^2} := \De^2$.
Integrating and making the change of variables $u=\nu - r\mu$, we obtain, by
definition of $c_1$, $-\int_{\ga}^{\nu - r\mu} \frac{du}{\sqrt{u}} =
\int_{\al}^{r} \De(s) \ ds \Leftrightarrow 2\sqrt{\ga} - 2\sqrt{\nu - r\mu} =
\int_{\al}^{r}\De$. Then: $$\mu = \frac{1}{r}\Big[ \nu - \big( \sqrt{\ga} -
\frac{1}{2}\int_{\al}^{r}\De \big)^2 \Big] $$ (in the interval
$[\widetilde{\al},\widetilde{\be}]$, where $\frac{1}{r}$ makes sense). Because
$\psi$ satisfies $\psi|_{[0,\al]}\equiv 0$, $\mu|_{[0,\al]}\equiv 0$ holds. 

It remains only $\mu|_{[\be,1]}\equiv 1$, for some $\psi$ carefully choosen.
However, 
$$\mu=1 \iff (*) \ r = \nu - \big( \sqrt{\ga} - \frac{1}{2} 
\int_{\al}^{r}\De
\big)^2,\textrm{ for all }r\geq\be.$$ 
By definition, $\De = \frac{1}{\sqrt{1-r}}$, for
$r\geq\be$. In particular, if $r\geq\be$, $\int_{\al}^r\De = \int_{\al}^{\be}\De
+ \int_{\be}^r\De = \int_{\al}^{\be}\De + 2\sqrt{1-\be} - 2\sqrt{1-r}$. Then
$(*)\iff \frac{1}{2}\int_{\al}^{\be}\De = \sqrt{\ga}-\sqrt{1-\be}$. 

Fix
$\ep>0$ small and consider $\psi_1$ s.t. 
$\psi_1|_{[\al+\ep,\be]}\geq 1-\ep$,
$\psi_1$ strictly increasing in $[\al,\be]$ and $\psi_1|_{[0,\al]}\equiv 0$,
$\psi_1|_{[\be,1]}\equiv 1$.
For a sufficiently small $\ep>0$, if $\rho_1$ and
$\De_1$ denotes the functions
associated to $\psi_1$, then  
$\rho_1|_{[\al+\ep,\be]}\geq 1- 2\ep$. In particular, $$\De_1\geq (1-3\ep)
\frac{1}{\sqrt{1-r}} \Rightarrow\frac{1}{2}\int_{\al}^{\be}\De_1\geq
(1-4\ep)\cdot\{ \sqrt{1-\al}-\sqrt{1-\be} \}.$$
By hypothesis, $\ga<1-\al$ that
imply $\frac{1}{2}\int_{\al}^{\be}\De_1 \geq \sqrt{\ga}-\sqrt{1-\be}$. With a
similar argument, we can take $\psi_0$ s.t. $\psi_0|_{[0,\al]}\equiv 0$, 
$\psi_0|_{[\be,1]}\equiv 1$ and $\psi_0|_{[\al,\be-\ep]}\leq\ep$. If $\ep>0$ is
small, we get $\frac{1}{2}\int_{\al}^{\be}\De_0 \leq \sqrt{\ga}-\sqrt{1-\be}$,
where $\De_0$ is the associated function to $\psi_0$. Now, consider the linear
combination $\psi_t = (1-t)\psi_0 + t\psi_1$. An easy verification is that
$\psi_t$ has the desired values in the intervals $[0,\al]$ and $[\be,1]$, for
any $t\in[0,1]$. By a continuity argument, it is 
not difficult that there is some
$t_0$ such that $\frac{1}{2}\int_{\al}^{\be}\De_{t_0} = \sqrt{\ga}-\sqrt{1-\be}$.
Finally, we consider $\psi=\psi_{t_0}$.

The Gauss-Kronecker curvature has the
properties of the lemma for the previous $\psi$.
Indeed, by the definitions, $rank(\vfi)=k$, since
$H_n$ is positive everywhere except at the cylindrical piece $S^k\times
D^{n-k}$ and $S^k\times D^{n-k}$ has exactly $k$ positive principal curvatures.
Moreover, $\psi>0$ if $r>\al$ implies that $\{ H_n=0 \} = S^k\times D^{n-k}$.
At this point, the informations about the radii of $S^k$ and $D^{n-k}$ are clear. This concludes the proof.   
{\cqd}

\begin{figure}[h]\label{fig1}
\begin{center}
\includegraphics[width=5.0in]{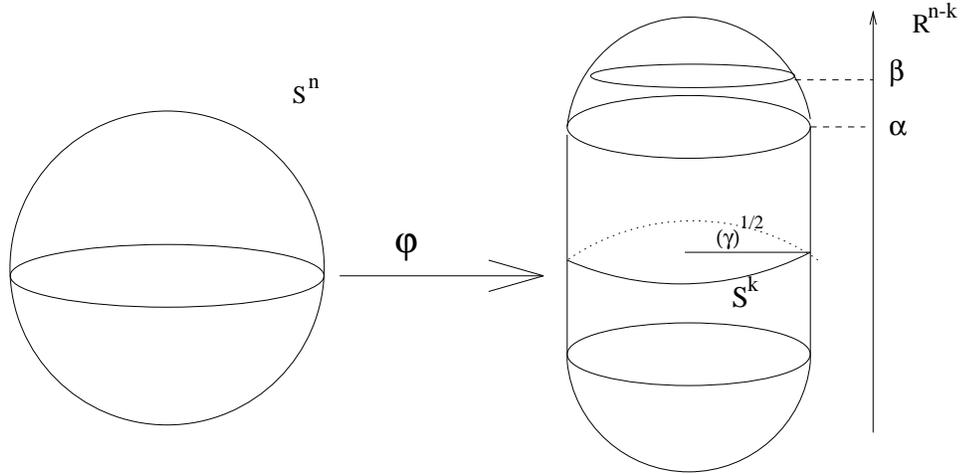}
\end{center}
\caption{Immersion with large region of zero curvature}
\end{figure}

The proof of lemma~\ref{im. cilindrica} is illustrated in figure 1.
The conditions on the numbers $\al,\be,\ga$ are necessary. In other case, our
immersions creates some undesired regions of negative Gauss-Kronecker. The role
of these numbers in the construction of $\vfi$ is explained in figure 1. 

In next section, we show how one can perturbate $\vfi$ to destroy regions of
zero curvature. With this is easy to obtain fractal immersions, by sucessively 
destroying ``large'' regions of zero curvature in such way that the set 
of zero curvature is fractal at the end of this process.    

\section{Removing zero Gauss-Kronecker curvature}
Before giving the statement and the proof of the main result of this paragraph,
we briefly describe the ideia of the lemma. Take the immersion $\vfi$ and an
open set $U$ of $S^{n-1}$. We want to perturbe $\vfi$ such that the region
correspondent to $U\times D^{1}$ has positive curvature (i.e., we want to 
``inflate'' 
the set), and this perturbation
should glue smoothly with the other region $(S^{n-1}-U)\times D^{1}$ of the 
cylinder $S^{n-1}\times D^{1}$. This can be done because the region of positive
curvature $\{ (z,w):|w|>\al\}$ permits the change of curvature in the ``goal''
region. Moreover, this method does not have points of negative curvature since
the curvature of $S^{n-1}$ is positive. The ideia is more clear contained in
figure 2 below.
 
\begin{figure}[h]\label{fig2}
\begin{center}
\includegraphics[width=5.0in]{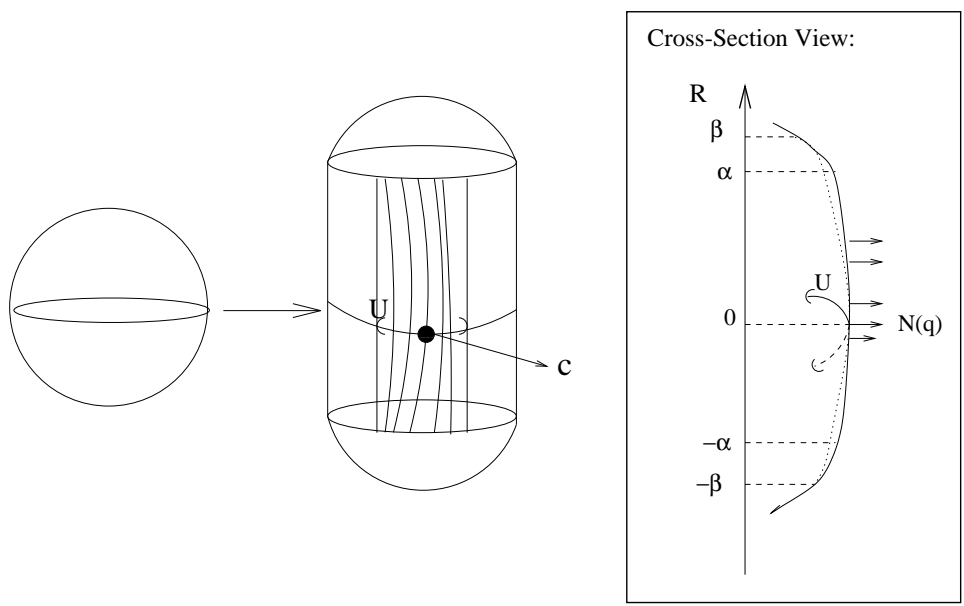}
\end{center}
\caption{Removing regions with zero curvature}
\end{figure}

With this figure in mind, the main lemma of this
section is :

\begin{lemma}\label{inflando telha} If $U$ is an open ball of $S^{n-1}$ with
radius bounded by $\de_0=\de_0(n-1)$, 
and $\vfi:S^n\rightarrow\real^{n+1}$ is the
immersion constructed above, then there exists a perturbation $\vfi_U$ of $\vfi$
such that $\vfi_U=\vfi$ in $(S^{n-1}-U)\times
D^{1} \ \bigcup \ \{(z,w): |w|\geq\be\}$ and the Gauss-Kronecker curvature of $\vfi_U(U\times D^{1})$ is
positive.
\end{lemma}

\proof[Proof of lemma~\ref{inflando telha}] Let $U\subset S^{n-1}$ be a round ball
of radius $\de$ and center $c$. Denote by
$d(z)=\frac{1}{2}\textrm{dist}(z,c)^2$, where $\textrm{dist}$ is the distance
function in $S^{n-1}$. We will consider a perturbation of 
$\vfi$ by normal variation, i.e., $F_t(q)=
q+t\cdot f(q)\cdot N(q)$, for $q\in\vfi(S^n)$ (here $N(q)$ denotes the normal
vector to $M:=\vfi(S^n)$ at $q$). To prove the result,
we need to show that a suitable $f:M\rightarrow \real$, $t>0$ have the desired
properties. In order to verify these properties, we first calculate the normal
vector $N_t(q_t)$ to $M_t:=F_t(M)$ at $q_t:=F_t(q)$. By definition, 

$$dF_t(q)v= v + t.f(q).dN(q)\cdot v + t.(df(q).v)\cdot N(q).$$  

In particular,
{\setlength\arraycolsep{2pt}
\begin{eqnarray}
0 & =& <\!\! N_t(q_t), dF_t(q).v\!\! >\nonumber\\
& = & <\!\! N_t(q_t),v\!\! > + t.(df(q).v).<\!\! N_t(q_t),N(q)\!\! > + \nonumber\\
& & + t.f(q).\!\! <\!\! N_t(q_t),dN(q).v\!\! >\nonumber
\end{eqnarray}}
$\forall \ v\in T_q M$. 

So, if $\de_0=\de_0(n-1)$ is small, we can define on $U$ a global ortonormal frame $\{ 
e_i(q) \}_{i=1}^{n}$ which diagonalizes $dN(q)$. If $k_i(q)$ are the respectives 
principal curvatures then the previous equation
implies:

$$<N_t(q_t),e_i(q)>=-\frac{t\cdot df(q)\cdot e_i(q)}{1+t.f(q).k_i(q)}\cdot 
<N_t(q_t),N(q)> .$$
 
Moreover, $||N_t(q_t)||^2=1$ says that:

$$<N_t(q_t),N(q)>=\Bigg[\sqrt{1+t^2\sum_{i=1}^{n}{\frac{(df(q)e_i(q))^2}
{(1+tf(q)k_i(q))^2}}}\Bigg]^{-1/2}:=\ze_t(q).$$

In particular, we can write:

$$N_t(q_t)=\ze_t(q).\Bigg[\sum_{i=1}^{n}-\frac{tdf(q)e_i(q)}{1+tf(q)k_i(q)}\cdot
e_i(q) + N(q)\Bigg] .$$

Differentiating the last expression:
{\setlength\arraycolsep{2pt}
\begin{eqnarray}
dN_t(q_t)\cdot v  & = & d\ze_t(q)\cdot v
\Big[\sum_{i=1}^{n}-\frac{tdf(q)e_i(q)}{1+tf(q)k_i(q)}\cdot e_i(q) + N(q) \Big] +\nonumber\\
& & +\ze_t(q)\cdot\Big[ \sum_{i=1}^{n}-\frac{tdf(q)e_i(q)}{1+tf(q)k_i(q)}\cdot
de_i(q)\cdot v + dN(q)\cdot v \Big] +\nonumber
\end{eqnarray}}

{\setlength\arraycolsep{2pt}
\begin{eqnarray}
& & +\ze_t(q)\cdot\Bigg[\sum_{i=1}^{n}-t\cdot\frac{e_i(q)}{(1+tf(q)k_i(q))}
\cdot\big[ <d^2f(q)\cdot v, e_i(q)> 
\nonumber\\
& & + <\textrm{grad} f(q),de_i(q)\cdot v> \big]\Bigg] - \nonumber\\
& & -\ze_t(q)\cdot\Bigg[\sum_{i=1}^{n}-t\cdot e_i(q)
\big[\frac{-tdf(q)e_i(q)}{(1+tf(q)k_i(q))^2}\cdot\big( df(q)v\cdot k_i(q)+\nonumber\\
& & +f(q)dk_i(q).v \big) \big] \Bigg]\nonumber
\end{eqnarray}}

Recall that we want to know the value of $\det dN_t(q_t)$. But it is not
difficult to see that, $\frac{\partial}{\partial t} \det dN_t(q_t)|_{t=0} =
(\frac{1}{c_2})^{n-1}\cdot <\!\!\frac{\partial}{\partial t}dN_t(q_t)|_{t=0}\cdot
e_n(q), e_n(q)\!\! >$ $+\sum_{i=1}^{n-1}(\frac{1}{c_2})^{n-2}\cdot \psi\cdot 
< \frac{\partial}{\partial t}dN_t(q_t)|_{t=0}\cdot
e_i(q), e_i(q)\!\! > $. In fact this follows from the fact that the determinant 
$\det A$ is the multilinear alternating $n$-form of the columns vectors 
$\det(A.e_1,\dots ,A.e_n)$. Thus,

$$\frac{d}{dt}\det A(t) |_{t=0} =\sum_{i=1}^{n}\det (A(t).e_1,\dots
,A'(t).e_i,\dots ,A(t).e_n)|_{t=0}.$$

Since $e_i(q)$, for $i=1,\dots ,n-1$, are eigenvectors of $dN(q)$ with eigenvalue
$\frac{1}{c_2}$ and $e_n$ is an eigenvector of $dN(q)$ with eigenvalue $\psi$
(see proof of lemma~\ref{im. cilindrica}), if we define $A(t)=
dN_t(q_t)$, $e_i=e_i(q)$, the formula above gives the claim.

We observe that the formula above uses (implicitly) that the codimension $(n-k)$ 
of the sphere $S^k$ is $1$. See remark~\ref{codim. 1} below. However,  
{\setlength\arraycolsep{2pt}
\begin{eqnarray}
\frac{\partial}{\partial t}
dN_t(q_t)|_{t=0}\cdot v & = & \frac{\partial}{\partial t}\ze_t(q)|_{t=0}\cdot
dN(q)\cdot v +\nonumber\\
& & +\lim\limits_{t\rightarrow 0}\ze_t(q)\cdot \Big[ -\sum_{i=1}^{n}
df(q)e_i(q) de_i(q)\cdot v \nonumber\\
& & -\sum_{i=1}^{n} e_i(q)\cdot(<\!\! d^2f(q)v,
e_i(q)\!\! >+<\!\! \textrm{grad}
f(q), de_i(q)v\!\! >)\Big] +\nonumber\\
& & +\frac{\partial}{\partial t}d\ze_t(q)|_{t=0}\cdot v\cdot N(q). 
\nonumber
\end{eqnarray}}

Taking $v=e_i(q), i=1,\dots , n$, we have :
$$<\!\!\frac{\partial}{\partial t} dN_t(q_t)|_{t=0}.e_i(q),e_i(q)\!\! >=-
<d^2f(q)e_i(q), e_i(q)>.$$

By the geometry of our immersion, we choose
$f(q)=f(z,w)=l_0. \la(d(z))\cdot\si(w)$ (see figure 2). 
Here $\si$ is a concave
function ($\si ''<0$) s.t. $\si(0)=0$, e.g., $\si(w)=-w^2/2$ (at least in $[0,\al_0]$, $\al<\al_0<\be$ close to $\al$) and $\la$ is
a bump function s.t. $\la\equiv 1$ if $t\leq 0$, $\la\equiv 0$ if $t\geq\de$,
$\la= e^{-1/\de -t}$ if $t$ is close to $\de$ and $\la$ is strictly decreasing
in $(0,\de )$.

An easy calculation shows that $$<d^2f(q)e_n(q),
e_n(q)>=l_0.\la(d(z))\cdot \si ''(w)\cdot <\frac{\partial}{\partial w},
e_n(q)>^2.$$
If $\al_0$ is sufficiently close to $\al$: $$<\frac{\partial}{\partial w},
e_n(q)>^2 \geq 1/2\textrm{ and } <d^2f(q)e_i(q), e_i(q)>=l_0.\si(w)\cdot\big[ \la''
||\textrm{grad} (d)||^2 + \la' \De d\big].$$ 

Then, 
$$\frac{\partial}{\partial t} \det dN_t|_{0}=l_0. 
(\frac{1}{c_2})^{n-2}.\big[ (\frac{1}{c_2})<\!\!\frac{\partial}{\partial w},
e_n(q)\!\!>^2.\la.\si ''+\psi\si\{\la ''||\textrm{grad}(d)||^2 +\la '\De d
\}\big].$$

To complete the proof, we show that the last expression is positive in $\{
(z,w): d(z)<\de, |w|\leq\al_0\}$ and it is small in $\{ (z,w): |w|\geq\be\}$.
This is sufficient because the derivative is positive imply that the curvature increases in the
``goal'' region, and the derivative is small (possible negative) does not creates
regions of negative curvature since the curvature starts positive in the
construction.

Since $\frac{1}{c_2}\geq \sqrt{\ga}$, $<\!\!\frac{\partial}{\partial w},
e_n(q)\!\!>^2\geq\frac{1}{2}$, $l_0$ is arbitrarily small and the term: $$
(\frac{1}{c_2})\cdot \la . \si '' . <\!\!\!\frac{\partial}{\partial w},
e_n(q)\!\!>^2 +\psi\si\{\la ''||\textrm{grad}(d)||^2 +\la '\De d \}$$
is uniformly bounded, the region $\{ (z,w): |w|\geq\be\}$ remains with positive
curvature. At the critcal region $\{
(z,w): d(z)<\de, |w|\leq\al_0\}$, we consider two cases:

\begin{enumerate}
\item  If $d(z)$ is close to $\de$, $\la '' = \big(\frac{1}{(\de -t)^4}- 2\cdot\frac{1}{(\de
-t)^3}\big)\cdot \la$ and $\la '= -\frac{1}{(\de -t)^2}\cdot\la$. So, if $d(z)$
is close to $\de$, i.e., $\de - d(z)$ is close to zero, the term  $\la
''||\textrm{grad}(d)||^2 +\la '\De d$ is positive. Since, $\si, \ \si ''$ is 
negative and $\psi$ is positive, this concludes the first case.
\item If $d(z)$ is far from $\de$, the term is $\la ''||\textrm{grad}(d)||^2 +
\la '\De d$ is bounded and $\la$ is positive and far away from zero. So, if
$\al_0$ is sufficiently close to $\al$, $\psi$ is small (and positive). This
concludes the second case.
\end{enumerate}
This completes the proof.
{\cqd}

\section{Proof of Theorem~\ref{imersao fractal}}

\proof[Proof of Theorem~\ref{imersao fractal}] Consider the good Cantor set 
$F=S^{n-1}-\bigcup\limits_{i=1}^{\infty} U_i$, where $U_i$ are round balls. We
fix an immersion $\vfi=\vfi_0$ given by lemma~\ref{im. cilindrica}. For each
$i$, let
$\vfi_i=\vfi_{U_i}$ a perturbation of $\vfi$ with support $\{ (z,w): d(z)<\de_i,
|w|\leq\al_0^i\}$, where $\de_i$ is the radius of $U_i$, $\al_0^i\geq\al$. Althought the existence
of $\vfi_i$ with the previous properties is not explicitly stated in
lemma~\ref{inflando telha}, this is contained in the proof. Finally, define
$x=\lim_n \vfi_n\circ\dots\circ\vfi_1\circ\vfi_0$. Observe that $U_i$ are pairwise
disjoint implies that the support of $\vfi_i$ are disjoint. So the limit
immersion $x$ above
exists and satisfies the desired properties.
{\cqd}

\section{Final Remarks}

We finish the paper with three remarks. The first remark is a possible 
generalization of
theorem A for open sets more general than round balls. The second remark is an
explanation about the restrictive hypothesis on the codimension. Finally, the
third is a generalization of the examples for the non-compact case. 

\begin{remark}We can replace in theorem~\ref{imersao fractal} the round balls by
sets with ``distance functions'' (i.e., bump functions with support \emph{equal}
to the open set) whose gradient is positive and bounded Laplacian. This follows 
from a carefull read of lemma~\ref{inflando telha}, 
the only place where properties of distance functions were used. 
\end{remark}

\begin{remark}\label{codim. 1}The proof of lemma $3.1$ works for codimension $1$ since for
higher codimensions, the determinat formula has a critical point of order
$(n-k-1)$ (the determinant is morally $ t^{n-k-1}\cdot\det A_t$, where $A_t$
are positive matrices close to $A_0$). So, our trick of calculate the first 
derivative
of this family and shows that the determinant increases does not work (the
critical point is ``flat'').
\end{remark}

\begin{remark} We observe that Matheus-Oliveira have also theorems for
non-compact manifolds. In this case, our techniques can be applied to give
examples of non-compact manifolds (or equivalently, immersions) satisfying
Matheus-Oliveira but not do Carmo-Elbert and Barbosa-Fukuoka-Mercuri hypothesis.
Since the proof is essentially analogous, we present here only a brief sketch of
this contruction :

\begin{itemize}

\item Consider the higher dimensional catenoid $M\subset\real^{n+1}$ (which is, topologically, a
sphere minus two points) and numbers $0\leq\al\leq\be$. As in lemma~\ref{im.
cilindrica}, one can modifythe catenoid to get a cylinder in the region $\{
(z,w) : |w|\leq\al \}$ and the original catenoid in $\{ (z,w):|w|\geq\be\}$
($(z,w)\in\real^{n}\times\real$). See
figure $3$ below :

\begin{figure}[h]\label{fig3}
\begin{center}
\includegraphics[width=5.0in]{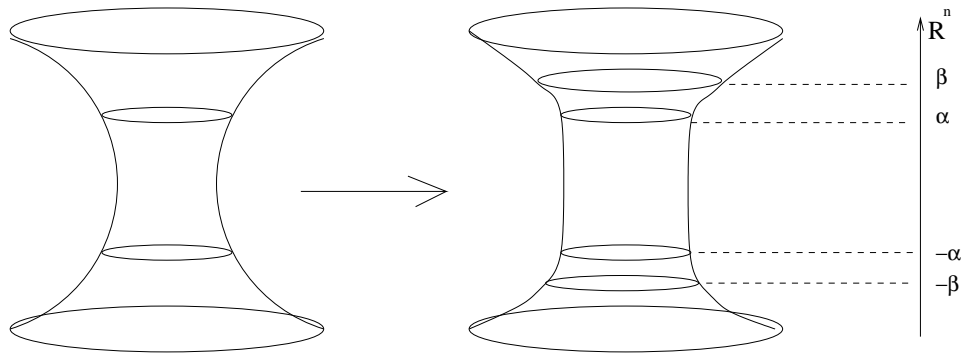}
\end{center}
\caption{The non-compact case}
\end{figure}

\item Using a lemma similar to lemma~\ref{inflando telha}, one can perturb the
new ``catenoid'' to obtain a fractal set of points of \emph{negative}
Gauss-Kronecker curvature (although the lemma~\ref{inflando telha} gives
positive curvature, we can obtain negative curvature if we choose in the proof
of lemma $3.1$, $\si$ a \emph{convex} function, because the catenoid has
negative curvature).

\end{itemize}

We point out that, in~\cite{MO}, there exist a statement about minimal
immersions but our method is \textbf{not} able to exhibit examples of minimal
hypersurfaces in Matheus-Oliveira's hypothesis but not in
Barbosa-Fukuoka-Mercuri's hypothesis.

\end{remark}

\textbf{Acknowledgements.} We are thankful to M. Viana for 
helpful advices and suggestions. 
Also to M. do Carmo and W. Santos for their interest in explicit
formulas for 
our ``geometrical description'' of fractal immersions during the geometry
seminar at IMPA. In addition, we are grateful to IMPA and his staff.     


\vspace{1cm}

\noindent    Alexander Arbieto ( alexande{\@@}impa.br )\\
             Carlos Matheus ( matheus{\@@}impa.br )\\
             IMPA, Est. D. Castorina 110, Jardim Bot\^anico, 22460-320 \\
             Rio de Janeiro, RJ, Brazil

\end{document}